\documentclass[12pt]{amsart}
\def\tS{\tilde{S}}
\def\N{{\Bbb N}}
\def\Z{{\Bbb Z}}

\newtheorem{Theorem}{Theorem}[section]
\newtheorem{Corollary}[Theorem]{Corollary}
\newtheorem{Lemma}[Theorem]{Lemma}

\theoremstyle{definition}

\theoremstyle{remark}

\begin{document}
\sloppy
\title{On the center of a Coxeter group}
\author{Tetsuya Hosaka} 
\address{Department of Mathematics, Utsunomiya University, 
Utsunomiya, 321-8505, Japan}
\date{September 9, 2005}
\email{hosaka@cc.utsunomiya-u.ac.jp}
\thanks{
Partly supported by the Grant-in-Aid for Young Scientists (B), 
The Ministry of Education, Culture, Sports, Science and Technology, Japan.
(No.~15740029).}
\maketitle
\begin{abstract}
In this paper, 
we show that the center of every Coxeter group is finite and 
isomorphic to $(\Z_2)^n$ for some $n\ge 0$.
Moreover, 
for a Coxeter system $(W,S)$, 
we prove that $Z(W)=Z(W_{S\setminus\tilde{S}})$ and $Z(W_{\tilde{S}})=1$, 
where $Z(W)$ is the center of the Coxeter group $W$ 
and $\tilde{S}$ is the subset of $S$ such that 
the parabolic subgroup $W_{\tilde{S}}$ is 
the {\it essential parabolic subgroup} of $(W,S)$ 
(i.e.\ $W_{\tilde{S}}$ is the minimum parabolic subgroup 
of finite index in $(W,S)$).
The finiteness of the center of a Coxeter group 
implies that a splitting theorem holds for Coxeter groups.
\end{abstract}

\section{Introduction and preliminaries}

In this paper, we investigate the center of a Coxeter group.
A {\it Coxeter group} is a group $W$ having a presentation
$$\langle \,S \, | \, (st)^{m(s,t)}=1 \ \text{for}\ s,t \in S \,
\rangle,$$ 
where $S$ is a finite set and 
$m:S \times S \rightarrow \N \cup \{\infty\}$ is a function 
satisfying the following conditions:
\begin{enumerate}
\item[(1)] $m(s,t)=m(t,s)$ for each $s,t \in S$,
\item[(2)] $m(s,s)=1$ for each $s \in S$, and
\item[(3)] $m(s,t) \ge 2$ for each $s,t \in S$
such that $s\neq t$.
\end{enumerate}
The pair $(W,S)$ is called a {\it Coxeter system}.
Let $(W,S)$ be a Coxeter system.
For a subset $T \subset S$, 
$W_T$ is defined as the subgroup of $W$ generated by $T$, 
and called a {\it parabolic subgroup}.
It is known that 
$(W_T,T)$ is also a Coxeter system (cf.\ \cite{Bo} and \cite{Hu}).
A subset $T\subset S$ is called a {\it spherical subset} of $S$, 
if the parabolic subgroup $W_T$ is finite.

The purpose of this paper is to prove the following theorems.

\begin{Theorem}\label{Thm1}
The center of every Coxeter group is finite and 
isomorphic to $(\Z_2)^n$ for some $n\ge 0$.
\end{Theorem}

A Coxeter system $(W,S)$ is said to be {\it irreducible}, 
if for any nonempty and proper subset $T$ of $S$, 
$W$ does not decompose as 
the direct product of $W_T$ and $W_{S \setminus T}$.

Let $(W,S)$ be a Coxeter system. 
Then there exists a unique decomposition $\{S_1,\ldots,S_r\}$ of $S$ 
such that $W$ is the direct product of 
the parabolic subgroups $W_{S_1},\ldots,W_{S_r}$ and 
each Coxeter system $(W_{S_i},S_i)$ is irreducible (cf.\ \cite{Bo} and \cite{Hu}).
We define $\tS=\bigcup \{S_i \,|\, W_{S_i} \ \text{is infinite} \}$, 
and the parabolic subgroup $W_{\tS}$ is called 
the {\it essential parabolic subgroup} of $(W,S)$.
We note that $W=W_{\tS}\times W_{S\setminus \tS}$ and 
$W_{S\setminus \tS}$ is finite.
In \cite{Ho2}, 
it was proved that 
the essential parabolic subgroup $W_{\tS}$ 
is the minimum parabolic subgroup of finite index in $(W,S)$.

We denote the center of a group $G$ by $Z(G)$.

We also prove the following theorem in Section~2.

\begin{Theorem}\label{Thm2}
For a Coxeter system $(W,S)$, 
$Z(W)=Z(W_{S\setminus\tilde{S}})$ and $Z(W_{\tilde{S}})=1$.
\end{Theorem}

For an irreducible Coxeter system $(W,S)$, 
if $W$ is infinite, then $W=W_{\tilde{S}}$.
Hence Theorem~\ref{Thm2} implies the following.

\begin{Corollary}\label{Cor1}
For an irreducible Coxeter system $(W,S)$, 
if the Coxeter group $W$ is infinite, 
then the center of $W$ is trivial.
\end{Corollary}

In \cite{Ho4}, 
we have obtained some splitting theorems 
for CAT(0) groups whose centers are finite.
Theorem~\ref{Thm1} and \cite[Theorem~2]{Ho4} 
implies the following splitting theorem for Coxeter groups.

\begin{Corollary}\label{Cor2}
Let $(W,S)$ be a Coxeter system and let $W=W_{S_1}\times W_{S_2}$.
Suppose that the Coxeter group $W$ acts geometrically 
on a CAT(0) space $X$. 
Then 
there exists a closed, convex, $W$-invariant, 
quasi-dense subspace $X'\subset X$ such that 
$X'$ splits as a product $X_1 \times X_2$ and 
the action of $W=W_{S_1}\times W_{S_2}$ on $X'=X_1 \times X_2$ 
is the product action.
\end{Corollary}

By \cite[Lemma~II.6.24]{BH}, we also can obtain the following corollary.

\begin{Corollary}\label{Cor3}
Suppose that a Coxeter group $W=W_{S_1}\times W_{S_2}$ 
acts geometrically on a CAT(0) space $X$.
Then $W_{S_1}$ and $W_{S_2}$ are convex-cocompact.
\end{Corollary}

Here the definition and some properties of 
{\it convex-cocompactness} is found in \cite{Ho1} and \cite{Ho3}.
We note that ``geometrically finiteness'' in \cite{Ho1} 
and ``convex-cocompactness'' in \cite{Ho3} coincide.

\section{Proof of the main theorems}

Let $(W,S)$ be a Coxeter system and $w\in W$.
A representation $w=s_1\cdots s_l$ ($s_i \in S$) is said to be 
{\it reduced}, if $\ell(w)=l$, 
where $\ell(w)$ is the minimum length of 
word in $S$ which represents $w$.

The following lemmas are known.

\begin{Lemma}[\cite{Bo}]\label{lem1}
Let $(W,S)$ be a Coxeter system. Suppose that $W$ is finite.
Then there exists a unique element $w_0\in W$ of longest length, 
and for each $w\in W$, $\ell(w_0w)=\ell(w_0)-\ell(w)$.
In particular, $w_0^2=1$.
\end{Lemma}

\begin{Lemma}[{\cite{Bo}, \cite[Lemma~7.11]{D1}}]\label{lem2}
Let $(W,S)$ be a Coxeter system, 
let $T\subset S$ and let $w\in W_T$. 
Then the following statements are equivalent.
\begin{enumerate}
\item[(1)] $W_T$ is finite and 
$w$ is the element of longest length in $W_T$.
\item[(2)] $\ell(wt)<\ell(w)$ for each $t\in T$.
\end{enumerate}
\end{Lemma}

Using lemmas above, 
we prove the following main theorem.

\begin{Theorem}
Let $(W,S)$ be a Coxeter system and 
let $Z(W)$ be the center of $W$.
\begin{enumerate}
\item[(1)] For each $w\in Z(W)$, 
there exists a spherical subset $T$ of $S$ such that 
$w$ is the element of longest length in $W_T$.
\item[(2)] $w^2=1$ for any $w\in Z(W)$.
\item[(3)] $Z(W)$ is finite.
\item[(4)] $Z(W)$ is isomorphic to $(\Z_2)^n$ for some $n\ge 0$.
\item[(5)] $Z(W)=Z(W_{S\setminus\tilde{S}})$.
\item[(6)] $Z(W_{\tilde{S}})$ is trivial.
\end{enumerate}
\end{Theorem}

\begin{proof}
Let $(W,S)$ be a Coxeter system and let $Z(W)$ be the center of $W$.

(1) Let $w\in Z(W)$ and 
let $w=s_1\cdots s_l$ be a reduced representation.
Then $ws_1=s_1w$, since $w\in Z(W)$.
Hence 
$$ ws_1=s_1w=s_1(s_1s_2\cdots s_l)=s_2\cdots s_l.$$
Thus $\ell(ws_1)<\ell(w)$ and 
$w=(s_2\cdots s_l)s_1$ is reduced.
Since $w\in Z(W)$, $ws_2=s_2w$.
Hence 
$$ ws_2=s_2w=s_2((s_2s_3\cdots s_l)s_1)=(s_3\cdots s_l)s_1.$$
Thus $\ell(ws_2)<\ell(w)$ and 
$w=(s_3\cdots s_l)s_1s_2$ is reduced.
By iterating the above argument, 
we obtain that 
$\ell(ws_i)<\ell(w)$ for each $i=1,\dots,l$.
Let $T=\{s_1,\dots,s_l\}$.
By Lemma~\ref{lem2}, 
$W_T$ is finite and 
$w$ is the element of longest length in $W_T$.

(2) By (1) and Lemma~\ref{lem1}, 
we have that $w^2=1$ for any $w\in Z(W)$.

(3) For a spherical subset $T$, 
let $w_T$ be the element of longest length in $W_T$.
By (1), 
$$Z(W)\subset \{w_T\,|\, \text{$T$ is a spherical subset of $S$}\}, $$
which is finite, since $S$ is finite.
Hence the center $Z(W)$ is finite.

(4) We note that $w^2=1$ for any $w\in Z(W)$ by (2) and 
$vw=wv$ for any $v,w\in Z(W)$ because $Z(W)$ is the center.
Thus $Z(W)$ is isomorphic to $(\Z_2)^n$ for some $n\ge 0$.

(5) Let $w\in Z(W)$, let $w=s_1\cdots s_l$ be a reduced representation
and let $T=\{s_1,\dots,s_l\}$.
Then $w$ is the element of longest length in $W_T$ by (1).
Let $s\in S\setminus T$.
Then $sw=ws$, since $w\in Z(W)$.
Hence $sws=w$ and 
$s(s_1\cdots s_l)s=s_1\cdots s_l$.
Here we note that $s\not\in T=\{s_1,\cdots,s_l\}$.
By Tits's theorem in \cite{T} and \cite[p.50]{Br}, 
$ss_i=s_is$ for each $i=1,\dots,l$.
This means that $st=ts$ for any $t\in T$ and $s\in S\setminus T$.
Hence $W$ splits as the product $W=W_T\times W_{S\setminus T}$.
Since $W_T$ is finite, $T\subset S\setminus \tS$ 
by the definition of $\tS$.
Hence $w\in W_T\subset W_{S\setminus \tS}$ for each $w\in Z(W)$.
Thus $Z(W)\subset W_{S\setminus \tS}$.
We note that $W=W_{\tS}\times W_{S\setminus \tS}$ and 
$Z(W)=Z(W_{\tS})\times Z(W_{S\setminus \tS})$.
Therefore $Z(W)=Z(W_{S\setminus \tS})$.

(6) We can obtain that $Z(W_{\tilde{S}})$ is trivial from (5).
\end{proof}

%

%
\end{document}